\documentclass[11pt]{amsart}

\newtheorem{theorem}{Theorem}[section]
\newtheorem{proposition}{Proposition}[section]
\newtheorem{lemma}{Lemma}[section]

\theoremstyle{definition}

\numberwithin{equation}{section}

\newcommand{\Real}{{\mathbb R}}
\newcommand{\Realm}{{\mathbb{R}^{\times}}}
\newcommand{\Com}{{\mathbb C}}
\newcommand{\G}{{\mathbf G}}
\newcommand{\half}{{\textnormal{\tiny{$\frac{1}{2}$}}}}

\begin{document}
\title[M\"{o}bius-convolutions and R.H.]{M\"{o}bius-convolutions and the Riemann hypothesis}
\author{Luis B\'{a}ez-Duarte}
\date{1:08 am, 22 November 2004; comments added on 19 April 2005.}
\maketitle
\begin{center}{22 November 2004}\end{center}

\begin{abstract}
The well-known necessary and sufficient criteria for the Riemann hypothesis of M. Riesz and of Hardy-Littlewood, based on the order of certain entire functions on the positive real axis, are here imbedded in a general theorem for a class of entire functions, which in turn is seen to be a consequence of a rather transparent convolution criterion. Some properties of the convolutions involved sharpen what is hitherto known for the Riesz function.
\end{abstract}

\section{Introduction}
RH stands for the Riemann hypothesis, and RHS for RH and simple zeros.  In this note the expression $f(x)\ll x^{a+\epsilon}$ always means
for $f(x)\ll_\epsilon x^{a+\epsilon}$ as $x\rightarrow+\infty$ for all $\epsilon>0$.

M. Riesz \cite{riesz} proved the following criterion:
\begin{theorem}\label{riesztheorem}
If
\begin{equation}\label{rieszfunction}
R(x):=\sum_{n=1}^\infty\frac{(-1)^{n+1}x^n}{(n-1)!\zeta(2n)},
\end{equation}
then
\begin{equation}\label{rieszcriterion}
RH\Longleftrightarrow R(x) \ll x^{\frac{1}{4}+\epsilon}.
\end{equation}
\end{theorem}

Likewise, G. H. Hardy and J. E. Littlewood \cite{hardy}, modifying Riesz's original idea, proved a similar theorem:

\begin{theorem}\label{hardytheorem}
\begin{equation}\label{hardyfunction}
H(x):=\sum_{n=1}^\infty \frac{(-x)^n}{n!\zeta(2n+1)},
\end{equation}
then 
\begin{equation}\label{hardycriterion}
RH \Longleftrightarrow H(x) \ll  x^{-\frac{1}{4}+\epsilon}.
\end{equation}
\end{theorem}
We shall see easily below that both theorems are included in a general RH criterion for a class of entire functions, which in turn is derived from a rather transparent, not to say trivial criterion for RH, Theorem \ref{grhthm}, predicated on the order at infinity of convolutions with $g(x):=\sum_{n\leq x}\mu(n)n^{-1}$. As a consequence, the known properties of both $R(x)$ and $H(x)$ can be sharpened. For example it is shown that both have an infinite number of positive, real zeros, in contrast with Riesz's statement that $R(x)$ has at least one such zero.

While approaching these theorems the author is well aware, not without trepidation, of E. C. Titchmarsh's trenchant comment in \cite{titch} and \cite{titchheath}:
\ \\

``\textit{These conditions have a superficial attractiveness since they depend explicitly only on values taken by $\zeta(s)$ at points in $\sigma>1$; but actually no use has ever been made of them.}"
\ \\

It has also been stated by J. M. Borwein, D. M. Bradley, and R. E. Crandall in \cite{borwein}, in reference to both Riesz's and Hardy-Littlewood's criteria, that 
\ \\

``\textit{It is unclear whether there be any computational value whatsoever to these equivalencies, especially as the big-O statement is involved and therefore infinite computational complexity is implicit, at least on the face of it. Still if there be any reason to evaluate such sums numerically, the aforementioned methods for recycling of $\zeta(\textnormal{even})$ or $\zeta(\textnormal{odd})$ values would come into play. }"
\ \\

In this context we have recently proved a discrete version of Riesz criterion in \cite{baez3} that lends itself well to calculations\footnote{K. Maslanka (personal communication) has carried out for us extensive numerical work that fits well with RH being true.}. It should seem appropriate to provide a general recipe to transform the power series criteria treated below to convert them to sequential criteria likewise.

\section{Preliminairies and notation}
Throughout $s$ stands for a complex variable with $\sigma=\Re(s)$. If $J$ is an interval, open, closed, or semi-closed, then $A J$ is the family of functions which are analytic in the strip $\sigma\in J$, and $A_c J$ is the family of  functions which are continuous in the strip and analytic in its interior.

In this note the \emph{(left)-Mellin transform} of $f:(0,\infty)\rightarrow\Com$ is defined by
\begin{equation}
f^\wedge(s):=\int_0^\infty t^{-s-1}f(t)dt,
\end{equation}
when the integral converges absolutely. For $a\in\Real$ introduce the associated norms\footnote{$N_a( . )$ is indeed the norm of $L_1((0,\infty),t^{-a-1}dt)$},  
$$
N_a(f):=\int_0^\infty t^{-a-1}|f(t)|dt.
$$
Clearly 
$$
N_0(\phi)=\|\phi\|_{L_1(\Realm)},
$$
where $\Realm:=(0,\infty)^\times$ is the multiplicative group of positive reals provided with its Haar measure $x^{-1}dx$. 

The following Lemma is standard.

\begin{lemma}\label{mellinanalytic}
If $N_{\sigma}(f)<\infty$ for $\sigma\in J$, then $f^\wedge\in A_c J$.
\end{lemma}

We shall say that $\phi$ is \textit{proper} when $N_\sigma(\phi)<\infty$ at least for $\sigma\in(-\frac{1}{2},1]$. In this case obviously $\phi^\wedge\in A_c(-\frac{1}{2},1]$. We say that $\phi$ is \textit{mellin-proper} if it is proper and $\phi^\wedge(s)\not=0$ in the strip $\sigma\in(-\frac{1}{2},0)$.

The \textit{Fourier transform} $\mathcal{F}$ in $f\in L_1(\Realm)$ is a continuous function which coincides with the Mellin transform on the line $\sigma=0$:
$$
\mathcal{F}(f)(\tau)=f^\wedge(i\tau), \ \ \tau\in\Real, \ \ f\in L_1(\Realm).
$$ 
\textbf{Note}: we reserve $\| . \|_p$ to denote the norm of the spaces $L_p((0,\infty),dx)$, except, of course, for $p=\infty$, which is the same for the measures $dx$ and $x^{-1}dx$. 

Let, as usual,
$$
M(x):=\sum_{n\leq x}\mu(n), \ g(x):=\sum_{n\leq x}\frac{\mu(n)}{n},
$$
and define
$$
g_1(x):=g(x^{-1})x^{-1}.
$$
Partial summation readily gives
\begin{equation}\label{growthMandg}
M(x)\ll x^{\alpha} \Longleftrightarrow g(x)
\ll x^{\alpha-1}, \ \ \alpha\in[\half,1),
\end{equation}
hence Littlewood's RH criterion becomes:
\begin{theorem}\label{grhthm}
\begin{equation}\label{grheq}
\fbox{$\displaystyle
RH \Longleftrightarrow g(x)\ll x^{-\frac{1}{2}+\epsilon},
$}
\end{equation}
and
\begin{equation}\label{grhs}
\fbox{$\displaystyle
g(x)\ll x^{-\frac{1}{2}} \Longrightarrow RHS.
$}
\end{equation}
\end{theorem}

The prime number theorem in the form $M(x)\ll x(\log x)^{-3}$ yields 
\begin{equation}\label{growthg}
g(x)\ll(\log x)^{-2},
\end{equation} 
hence
\begin{equation}\label{gL1dx/x}
g\in L_p(\Realm), \ \ p\in[1,\infty].
\end{equation}

Summing by parts the Dirichlet series of $(\zeta(s+1))^{-1}$ one obtains
\begin{equation}\label{gmellin}
\frac{1}{s\zeta(s+1)}=g^\wedge(s)=\int_0^\infty x^{-s-1}g(x)dx, \  \sigma\geq0,
\end{equation}
where the integral is absolutely convergent, and $g^\wedge\in A[0,\infty)$. Setting $s=0$ above yields
\begin{equation}\label{intg(x)x^{-1}}
\int_0^\infty x^{-1}g(x)dx=1.
\end{equation}
The Littlewood criterion (\ref{grheq}) for RH can be translated to a criterion based on $L_p$-norms, which may have some theoretical significance:

\begin{theorem}\label{grhlpthm}
With $g_1(x):=x^{-1}g(x^{-1})$ we have
\begin{equation}\label{grhlpeq}
\fbox{$\displaystyle
RH \Longleftrightarrow (\|g_1\|_p<\infty, \ \forall p\in[1,2)),
$}
\end{equation}
and, unconditionally, $\|g_1\|_2=\infty$.
\end{theorem}

\begin{proof}
Note that $\|g\|_1<\infty$ is unconditionally true. The direct implication follows directly from Theorem \ref{grhthm}. For the converse assume $\|g_1\|_p<\infty$ for some $p\in(1,2)$, then an application of H\"{o}lder's inequality gives
$$
\int_0^\infty x^{-\sigma-1}|g(x)|dx=\int_0^1 x^\sigma|g_1(x)|dx\leq
\|g_1\|_p\left(\int_0^1 x^{\sigma q}dx\right)^{\frac{1}{q}}<\infty, 
$$
for $\sigma>\frac{1}{p}-1$. Since this ranges in $(-\frac{1}{2},0)$ we apply Lemma \ref{mellinanalytic} to the integral on the right-hand side of (\ref{gmellin}) to obtain the analytic extension of $(s\zeta(s+1))^{-1}$, and, hence, RH. 

Now let $f(x):=\int_0^x g(t)dt$, then 
$\int_1^\infty x^{-s-1}f(x)dx=(s(s-1)\zeta(s))^{-1}$, a calculation easily justified at least for $\sigma\geq1$. The right-hand side has poles on the line $\sigma=\frac{1}{2}$, hence, by the order Lemma 2.1 in \cite{baez1}, $f(x)\not=o(x^{\frac{1}{2}})$, and Lemma 2.3 in \cite{baez1} implies $\|g\|_2=\infty$; but $\|g\|_2=\|g_1\|_2$.
\end{proof}

\section{Convolution criterion for the Riemann hypothesis}

\subsection{The convolution operator $\G$}
For measurable $\phi:(0,\infty)\rightarrow\Com$, and any given $x>0$, we define 
\begin{equation}\label{gdef}
\fbox{$\displaystyle
\G\phi(x):=\int_0^\infty g(x t)\phi(t^{-1})t^{-1}dt,
$}
\end{equation}
provided the integral converges absolutely. A rather general condition for existence is this: $\G\phi(x)$ exists and is continuous for all $x>0$ when $\phi$ is bounded on \textit{every} interval $(\delta,\infty)$, $\delta>0$. So far the most interesting class of examples arises as follows: for any power series of type
\begin{equation}\label{phipower}
\phi(z):=\sum_{n=1}^\infty a_n z^n,
\end{equation}
write
\begin{equation}\label{phistarpower}
\phi^\star(z):=\sum_{n=1}^\infty \frac{a_n z^n}{n\zeta(n+1)}.
\end{equation}
One must note the trivial fact that there exists another series $\phi_1(z)$ such that $\phi(z)=\phi_1^\star(z)$. At first blush the following proposition seems only interesting for entire functions, and for $z=x>0$; it is nevertheless convenient not to lose sight of the general case.
\begin{proposition}\label{analyticg}
If $\phi(z)$ and $\phi^\star(z)$ are as above, and $R>0$ is their common radius of convergence, then
\begin{equation}\label{phistarisgphi1}
\fbox{$\displaystyle
\phi^\star(z)=\int_0^1g_1(t)\phi(zt)dt, \ |z|<R,
$}
\end{equation}
and
\begin{equation}\label{phistarisgphi2}
\fbox{$\displaystyle
\phi^\star(x)=\G\phi(x), \ 0<x<R.
$}
\end{equation}
\end{proposition}
\begin{proof}
Write (\ref{gmellin}) for $s=n$ as 
$(n\zeta(n+1))^{-1}=\int_0^1 g(t^{-1})t^{n-1}dt$, and substitute in the definition of $\phi^\star$ to obtain (\ref{phistarisgphi1}), since the interchange of sum and integral is totally trivial. Now take $z=x\in(0,R)$ and change variables to obtain
$$
\phi^\star(x)=\int_0^\infty g(xt)\phi(t^{-1})t^{-1}dt=\G\phi(x), \ 0<x<R.
$$
\end{proof}

Clearly the Riesz and Hardy-Littlewood functions defined in (\ref{rieszfunction}), and (\ref{hardyfunction}), satisfy
\begin{equation}\label{rieszG}
x^{-1}R(x^2)=\G\alpha(x), \ \ \ \alpha(x)=x(1-2x^2)e^{-x^2},
\end{equation}
and
\begin{equation}\label{hardyG}
H(x^2)=\G\beta(x), \ \ \ \beta(x)=-2x^2 e^{-x^2},
\end{equation}
for all $x>0$. It is quite easy to see that both $\alpha$ and $\beta$ are mellin-proper.
\ \\

We now establish in some generality the main elementary properties of $\G$, in a context relevant to RH.

\begin{lemma}\label{GL1}
Assume $\phi\in L_1(\Realm)$, then
$$
\G\phi=g\ast\phi\in L_1(\Realm)
$$
in the sense of $\Realm$-convolutions. $G\phi$ is continuous and vanishes both at $0$ and at $\infty$. In particular $G\phi$ is \emph{bounded}.
\begin{equation}\label{normgphiL1m}
N_0(\G\phi)\leq N_0(g)N_0(\phi),
\end{equation}
The left-Mellin transform of $\G\phi$ exists al least on the line $\sigma=0$:
\begin{equation}\label{mellingphi1}
\fbox{$\displaystyle
\int_0^\infty x^{-s-1}\G\phi(t)dx=\frac{1}{s\zeta(s+1)}
\int_0^\infty t^{-s-1}\phi(t)dt, \ \sigma=0,
$}
\end{equation}
where both integrals converge absolutely. A fortiori
\begin{equation}\label{intgphi}
\int_0^\infty x^{-1}\G\phi(x)dx=\int_0^\infty x^{-1}\phi(x)dx.
\end{equation}
\end{lemma}

\begin{proof}
$g\in L_1(\Realm)$, hence the integral in (\ref{gdef}) is absolutely convergent for a.e.$\hspace{1mm} x$, and $\G\phi\in L_1(\Realm)$. Since $g$ is bounded, the continuity, as well as the vanishing at the ends follow from Lebesgue's dominated convergence theorem. The inequality (\ref{normgphiL1m}) is just ad hoc notation for the Banach algebra property. On account of (\ref{gmellin}) 
$\mathcal{F}(\G\phi)=\mathcal{F}(g\ast\phi)=\mathcal{F}(g)\mathcal{F}(\phi), \ \phi\in L_1(\Realm)$
translates into (\ref{mellingphi1}), and letting $s=0$ one obtains (\ref{intgphi}).  
\end{proof}

It is quite natural to seek conditions to extend the range of the identity (\ref{mellingphi1}) to obtain a sufficient condition for RH. It proves convenient to identify separately some of the simple steps of the process.

\begin{lemma}\label{gphiin01}
If $N_0(\phi)<\infty$, then 
$$
h(s):=\int_0^1 x^{-s-1}\G\phi(x)dx
$$ 
is in $A_c(-\infty,0]$. 
\end{lemma}

\begin{proof}
Note that the analyticity in the half-plane $\sigma\in(-\infty,0)$ is trivial since $\G\phi$ is bounded by Lemma \ref{GL1}. But we need continuity at the boundary. So we argue that in the interval of integration we have
$x^{-\sigma}\leq1$ for all $\sigma\leq0$, hence by (\ref{normgphiL1m}). 
$$
\int_0^1 x^{-\sigma-1}|\G\phi(x)|dx\leq N_0(\G\phi)\leq N_0(g)N_0(\phi)<\infty.
$$
Now apply Lemma \ref{mellinanalytic}.
\end{proof}

\begin{lemma}\label{extlemma}
If $\phi$ is proper, and $\G\phi(x)\ll x^{-\frac{1}{2}+\epsilon}$, then
\begin{equation}\label{mellingphi2}
(\G\phi)^\wedge(s)=\frac{1}{s\zeta(s+1)}\phi^\wedge(s), \ \ \sigma\in(-\half,0],
\end{equation}
where both sides are in $A_c(-\frac{1}{2},0]$. 
\end{lemma}

\begin{proof}
On account of (\ref{mellingphi1}) and Lemma \ref{mellinanalytic} all we are required to do is to show that $N_\sigma(\G\phi)<\infty$ for $\sigma\in(-\frac{1}{2},0]$, then invoke analytic continuation. Accordingly we split the $N_\sigma(\G\phi)$ integral at $x=1$: the interval $(0,1)$ is already taken care of by Lemma \ref{gphiin01}. Now, for the interval $(1,\infty)$ we have 
$$
\int_1^\infty x^{-\sigma-1}|\G\phi(x)|dx \ll  
\int_1^\infty x^{-1-((\sigma+\frac{1}{2})-\epsilon)}dx<\infty.
$$
when $0<\epsilon<\sigma+\frac{1}{2}$.
\end{proof}

We therefore have the following \textit{necessary and sufficient condition for} RH:

\begin{theorem}[RH convolution criterion]\label{grhthm}
Let $\phi$ be proper, then 
\begin{equation}\label{necessity}
\fbox{$\displaystyle
RH\Longrightarrow G\phi(x)\ll x^{-\frac{1}{2}+\epsilon},
$}
\end{equation}
and
\begin{equation}\label{phimellin0}
\fbox{$\displaystyle
G\phi(x)\ll x^{-\frac{1}{2}+\epsilon}\Longrightarrow (\zeta(s+1)=0 \Rightarrow \phi^\wedge(s)=0),
$}
\end{equation}
where only $\sigma>-\frac{1}{2}$ is considered. \emph{A fortiori}, if $\phi$ is mellin-proper then
\begin{equation}\label{grhequiv}
\fbox{$\displaystyle
RH \Longleftrightarrow \G\phi(x)\ll x^{-\frac{1}{2}+\epsilon}.
$}
\end{equation}
\end{theorem}

It is easily seen that one also has the following general equivalence:
\begin{equation}\label{grhequiv1}
\fbox{$\displaystyle
RH \Longleftrightarrow (\G\phi(x)\ll x^{-\frac{1}{2}+\epsilon}, \forall \textnormal{ proper } \phi)
$}
\end{equation}

\begin{proof}[Proof of Theorem \ref{grhthm}]
The necessity implication (\ref{necessity}) follows using the Littlewood criterion (\ref{grheq}), and a simple estimate of the integral (\ref{gdef}) defining $\G\phi(x)$. The implication (\ref{phimellin0}) is just read off from Lemma \ref{extlemma}.  After this the necessary and sufficient condition (\ref{grhequiv}) becomes trivial.
\end{proof}

As a corollary we get the generalization of Riesz's criterion to entire functions.

\begin{theorem}[Entire function RH criterion]\label{generalriesz}
Let $\phi$ be an entire function vanishing at zero as in (\ref{phipower}) and $\phi^\star$ the associated entire function defined in (\ref{phistarpower}). If $\phi$ is mellin-proper\footnote{If an entire function $\phi$ vanishes at zero, and $\phi(x)\ll x^{-a}$ for some $a<\frac{1}{2}$, then it is proper.}, then
\begin{equation}\label{entirecriterion}
\fbox{$\displaystyle
RH \Longleftrightarrow \phi^\star(x) \ll x^{-\frac{1}{2}+\epsilon}.
$}
\end{equation}
\end{theorem}
In view of  (\ref{rieszG}) and (\ref{hardyG}) the above criterion immediately proves the Riesz and the Hardy-Littlewood criteria in Theorems \ref{rieszcriterion} and \ref{hardycriterion}.

\section{Further properties of $\G\phi$}

The same argument used to prove the necessity implication in the main Theorem \ref{grhthm} gives: $g(x)\ll x^{-\frac{1}{2}}$ and $N_{-1/2}(\phi) < \infty$ imply $\G\phi(x)\ll x^{-\frac{1}{2}}$. In view of (\ref{grhs}) it is more interesting to show the next implication:

\begin{theorem}\label{rhs}
If $\phi$ is mellin-proper, and additionally $\phi^\wedge\in A_c[-\frac{1}{2},1]$ does not vanish on the line $\sigma=-\frac{1}{2}$, then
\begin{equation}\label{rhssuff}
\fbox{$\displaystyle
\G\phi(x)\ll x^{-\frac{1}{2}} \Longrightarrow RHS.
$}
\end{equation}
\end{theorem}

\begin{proof}
By the hypotheses and Lemmas \ref{gphiin01} and \ref{extlemma} we have 
$$
\frac{1}{s\zeta(s+1)}\phi^\wedge(s)=h(s)+\int_1^\infty x^{-s-1}\G\phi(x)dx,
 \ \sigma\in[-\half,0],
$$
where $h(s)\in A_c(-\infty,0]$. Now assume 
$s_0=-\frac{1}{2}+i\tau$ is a zero of $\zeta(s+1)$. We take $s=\sigma+i\tau$ and let $\sigma\downarrow-\frac{1}{2}$. While this happens $|h(s)|<a<\infty$, and $|\phi(s)|>b>0$; hence the main assumption that $\G\phi(x)\ll x^{-\frac{1}{2}}$ yields 
$$
\left|\frac{1}{s\zeta(s+1)}\right| \leq \frac{a}{b}+\frac{1}{b}\int_1^\infty x^{-\sigma -1}|\G\phi(x)|dx \ll \int_1^\infty x^{\sigma-\frac{3}{2}}dx \ll \frac{1}{\sigma -\frac{1}{2}},
$$
which shows $s_0$ is a simple zero.
\end{proof}

Theorem \ref{grhthm} admits a variant quite like the $L_p$-criterion of Theorem \ref{grhlpthm}. Conditions on $\phi$ need to be a bit stronger than in the main RH criterion Theorem \ref{grhthm}, namely, they are set as in the above sufficiency Theorem \ref{rhs} for RHS.

\begin{theorem}\label{lpvariantthm}
Let $\phi$ be mellin-proper, and define
$$
\psi(x):=x^{-1}\G\phi(x^{-1}),
$$
then
\begin{equation}\label{lpvarianteq}
\fbox{$\displaystyle
RH \Longleftrightarrow (\|\psi\|_p<\infty, \ \forall p\in [1,2)),
$}
\end{equation}
Furthermore, unconditionally on RH, if $\phi^\wedge(s)\in A[-\frac{1}{2},0)$ does not vanish on the line $\sigma=-\frac{1}{2}$, then  
\begin{equation}\label{notinL2}
\|\psi\|_2=\infty
\end{equation}
\end{theorem}

\begin{proof}
First note that $\|\psi\|_1=N_0(\G\phi)<\infty$ is unconditionally true by (\ref{normgphiL1m}). Now assume RH. Pick any $p\in(1,2)$ and write
$$
\|\psi\|_p^p=\int_0^\infty|G\phi(x)|^p x^{p-2}dx=\int_0^1 + \int_1^\infty.
$$
$\G\phi$ is bounded by Lemma \ref{GL1} so the first integral is finite. For the second note that $RH\Longrightarrow g(x)\ll x^{-\frac{1}{2}+\epsilon} \Longrightarrow \G\phi(x)\ll x^{-\frac{1}{2}+\epsilon}$, so the integral is finite when $p<\frac{2}{1+2\epsilon}$. Conversely, assume $\|\psi\|_p<\infty$ for an arbitrary $p\in(1,2)$. Write the fundamental identity (\ref{mellingphi1}) as
$$
(\G\phi)^\wedge(s)=\frac{1}{s\zeta(s+1)}\phi^\wedge(s)=h(s)+\int_0^1 t^s \psi(t) dt, \ \ \sigma=0,
$$
where $h(s)$ is as in Lemma \ref{gphiin01}, thus in $A_c(-\infty,0]$. On the other hand, H\"{o}lder's inequality shows that the last integral above converges absolutely for $\sigma>\frac{1}{p}-1$, thus $(\G\phi)^\wedge(s)$ is analytic in the strip $\sigma\in(\frac{1}{p}-1,0]$, and this implies that $\zeta(s)$ does not vanish in $\sigma>\frac{1}{p}$. This means RH is true.

Finally, under the additional hypothesis for $\phi$, we show $\|\psi\|_2=\infty$ unconditionally. Define $f(x):=\int_0^x \G\phi(t)dt$.  The calculation 
$$
\int_1^\infty x^{-s-1}f(x)dx=\frac{\phi^\wedge(s-1)}{(s-1)s\zeta(s)}
$$ 
is easily justified for $\sigma=1$, where by hypothesis $\phi^\wedge(s-1)$ is in $A[\frac{1}{2},1)$, and does not vanish on the critical line, thus the left-hand side transform has a meromorphic continuation that certainly has poles on the critical line. Thus order Lemma 2.1 in \cite{baez1} can be applied, so $f(x)\not=o(x^{-\frac{1}{2}})$, and then Lemma 2.3 in \cite{baez1} yields $\|\G\phi\|_2=\infty$; but $\|\G\phi\|_2=\|\psi\|_2$.
 \end{proof}
\ \\

With mild restrictions on $\phi$, but \emph{without assumption of RH}, the order of $\G\phi(x)$, like that of $g$, is limited between $o(1)$ and $O(x^{-\frac{1}{2}})$. For the \emph{maximal order} it is clear that
\begin{equation}\label{maxorderGphi}
\G\phi(x)=o(1), \ x\rightarrow+\infty,
\end{equation}
whenever $\phi\in L_1(\Realm)$ by \ref{GL1}. Using the error term of the prime number theorem one could be more specific about the order implied in $o(1)$. At present we do not think this is worthwhile. 

As for the \emph{minimal order}: if $\phi$ satisfies the conditions in Theorem \ref{grhthm} \emph{farther to the left}, that is, for some $\delta>0$ one has both $N_\sigma(\phi)<\infty$, and $\phi^\wedge(s)\not=0$ for $\sigma\in(-\frac{1}{2}-\delta,0]$, then
\begin{equation}\label{maxorderGphi1}
\G\phi(x)\not\ll x^{-1/2-\epsilon},
\end{equation}
as follows easily from the reasoning in the proof of Theorem \ref{grhthm}. But \textit{more is true} without strengthening the hypotheses on $\phi$:

\begin{proposition}\label{minimalorder}
If $\phi$ is mellin-proper, then
\begin{equation}\label{gphinotlittleo}
\G\phi(x)\not= o(x^{-\frac{1}{2}}).
\end{equation}
\end{proposition}

\begin{proof}
On account of Lemmas \ref{gphiin01} and \ref{extlemma} we see that the integral
\begin{equation}\label{laplacetransform}
f(s):=
\int_1^\infty x^{-s-1}\G\phi(x)dx = -h(s)+\frac{1}{s\zeta(s+1)}\phi^\wedge(s),
\end{equation}
where $h(s)\in A(-\infty,0]$, has a finite abscissa of convergence $\alpha\geq-\frac{1}{2}$. In this half-plane there are poles for $f(s)$, therefore, by the \emph{order Lemma} 2.1 in \cite{baez1}, we get $\G\phi(x)\not=o(x^{-\frac{1}{2}})$.
\end{proof}

Finally we deal with the oscillations of $\G\phi$. If $\phi\in L_1(\Realm)$ is \emph{real}, and $\int_0^\infty \phi(t)t^{-1}dt=0$, then $\int_0^\infty \G\phi(t)t^{-1}dt=0$ by (\ref{intgphi}), so, unless $\phi=0$,  \emph{there must be at least one change of sign in} $x>0$. This is a simpler argument than Riesz's in \cite{riesz} to establish the existence of at least one sign change for $R(x)$. But again, one can prove more, namely, that the infinite number of oscillations of $g(x)$ are transmitted to $\G\phi(x)$ \emph{unconditionally}:

\begin{theorem}\label{gphioscillations1}
If $\phi$ is real and mellin-proper, then there exists a $\beta\in[0,\frac{1}{2})$ such that 
\begin{eqnarray}\label{gphioscillations2}
\liminf_{x\rightarrow\infty}x^{\beta+\epsilon}\G\phi(x)&=&-\infty,\\\nonumber
\limsup_{x\rightarrow\infty}x^{\beta+\epsilon}\G\phi(x)&=&+\infty,
\end{eqnarray}
for all $\epsilon>0$.
\end{theorem}

\begin{proof}
Consider the same integral $f(s)$ defined in (\ref{laplacetransform}) with abscissa of convergence $\alpha\in(-\frac{1}{2},0]$. It is obvious that $s=\alpha$ cannot be a singularity of $f(s)$, therefore the \emph{oscillation Lemma} 2.2 in \cite{baez1}, based on Landau's theorem for Laplace tranforms (see \cite{widder}), yields the conclusion with $\beta=-\alpha$.
\end{proof}
\ \\

Obviously the above results apply to Riesz's and Hardy-Littlewood's functions through the relations $x^{1}R(x)=\G\alpha(x)$ and $H(x^2)=\G\beta(x)$. This clearly extends and sharpens some of the properties of $R(x)$ stated by M. Riesz in \cite{riesz}. In a forthcoming installment of this series we shall furthermore explore some interesting properties of this and other special power series. We also take up the important subject of inverting $G$ and its connection to the M\"{u}ntz and M\"{o}bius operators $P$ and $Q$ (see \cite{baez2}), as well as to \textit{Burnol's co-Poisson intertwining}; this last topic has a vast literature by this author behind it, but which I here only refer, for the moment, to another note in this series, namely \cite{baez4} .

\bibliographystyle{amsplain}

\ \\
\ \\
\noindent Luis B\'{a}ez-Duarte\\
Departamento de Matem\'{a}ticas\\
Instituto Venezolano de Investigaciones Cient\'{\i}ficas\\
Apartado 21827, Caracas 1020-A\\
Venezuela\\
\email{lbaezd@cantv.net}

 \end{document}